\documentclass[12pt]{article}
\usepackage{amsmath,amsthm,amsfonts,amssymb}
\usepackage[all]{xy}
\newtheorem{thm}{Theorem}[section]
 \newtheorem{lemma}[thm]{Lemma}
\newtheorem{cor}[thm]{Corollary}
\newtheorem{proposition}[thm]{Proposition}

\newtheorem{prop}[thm]{Proposition}

\newtheorem{sublem}[thm]{Sublemma}
\theoremstyle{definition}
\newtheorem{defn}[thm]{Definition}

\newtheorem{notn}[thm]{Notation}
\newtheorem{com}[thm]{Remark}

\newtheorem{remark}[thm]{Remark}
\theoremstyle{remark}
\DeclareMathOperator{\Fix}{Fix}

\DeclareMathOperator{\Int}{Int}

\DeclareMathOperator{\Out}{Out}

\DeclareMathOperator\closure{Cl}
\DeclareMathOperator\mlog{mlog}
\DeclareMathOperator\Fr{Fr}
\DeclareMathOperator\rk{rk}

\newcommand{\F}{{\cal F}}

\newcommand{\mcg}{\cal MCG}
\newcommand{\T}{\cal T}

\newcommand{\CI}{\cal FI}
\newcommand{\inv}{{-1}}
\newcommand\from\colon
\newcommand\abs[1]{\left| #1 \right|}

\title{The expansion factors of \\ an outer automorphism and its inverse}
\author{Michael Handel\thanks{Supported in part by NSF grant DMS0103435.}\ \ 
and Lee Mosher\thanks{Supported in part by NSF grant
DMS0103208.}}

\def\sl3z{SL(3, \mathbb Z)}

\def\eg{EG}
\def\neg{NEG}

\newcommand\subgroup{<}
\newcommand\BookOne{\cite{BFH:TitsOne}}
\newcommand\BH{\cite{BestvinaHandel:tt}}
\newcommand\A{{\mathcal A}}

\newcommand\union\cup
\renewcommand\L{{\mathcal L}}
\newcommand\composed\circ
\newcommand\suchthat{\bigm|}
\DeclareMathOperator\rank{rank}
\newcommand\restrict{\bigm|}

\begin{document}
\maketitle

\begin{abstract} 
A fully irreducible outer automorphism $\phi$ of the free group $F_n$ of rank~$n$ has an expansion factor which often differs from the expansion factor of $\phi^\inv$. Nevertheless, we prove that  the ratio between the logarithms of the expansion factors of $\phi$ and $\phi^\inv$ is bounded above by a constant depending only on the rank $n$. We also prove a more general theorem applying to an arbitrary outer automorphism of $F_n$ and its inverse and their two spectrums of expansion factors. 
\end{abstract}

\section{Introduction}  We denote the free group on $n$ letters by $F_n$, the outer
automorphism group of $F_n$ by $\Out(F_n)$ and the mapping class group of a closed
surface  by $\mcg$.  There has been a great deal of success in developing the analogy
between $\Out(F_n)$  and $\mcg$ in general and between Outer Space $X_n$ and the
Teichmuller space $\T$ in particular. The survey paper \cite{Vogtmann:OuterSpaceSurvey} is an excellent
reference.  One notable exception is the lack in $\Out(F_n)$ of an analogue for
Teichmuller theory.  This is the first of several papers  on this topic; see also the
companion paper \cite{HandelMosher:parageometric}.   Our focus is not on finding an analogue of
the Teichmuller metric but rather on finding an analogue of Teichmuller geodesics.

If $\gamma$ is a Teichmuller geodesic that is $\psi$-invariant for some $\psi \in
\mcg$, then $\psi$ is pseudo-Anosov and $\gamma$ serves as an axis for the action of
$\psi$ on $\T$: the minimum translation distance of $\psi$ is realized precisely on
$\gamma$. This translation distance is $\log \lambda$, where $\lambda$ is the expansion
factor of $\psi$. Symmetry of the Teichmuller metric implies that $\psi$ and
$\psi^\inv$  have the same axis and the same expansion factor.

A strictly analogous result does not hold in $\Out(F_n)$ no matter what metric is placed on
$X_n$.   An outer automorphism $\phi$ is {\it fully irreducible}, written $\phi\in \CI$, if
there are no free factors of $F_n$ that are invariant under an iterate of $\phi$.  Fully
irreducible $\phi$ correspond to pseudo-Anosov mapping classes. They are the most well behaved
and best understood elements of $\Out(F_n)$.   The action induced by $\phi \in\CI$ on the
boundary of $X_n$ has exactly two fixed points, one a source and the other a sink with all other
orbits limiting on these points \cite{LevittLustig:NorthSouth}.  Each $\phi \in \CI$ has a well
defined expansion factor that can be thought of as the exponential growth rate of the action of
$\phi$ on conjugacy classes in $F_n$.  In contrast to $\mcg$, the expansion factor $\lambda$ of
$\phi$ need not equal the expansion factor $\mu$ of $\phi^\inv$.  This implies that the action
of $\phi$ on $X_n$ cannot have an axis in the traditional sense: if there were a metric on
$X_n$ and an axis for $\phi$ that uniquely minimizes translation length, symmetry of the metric
would imply that $\lambda=\mu$ as in Teichmuller space.

Any attempt to develop a Teichmuller theory for $\Out(F_n)$ must address this phenomenon. 
Given that the expansion factor of $\phi^n$ is $n$ times the expansion factor of $\phi$, the
correct way to measure the variation between $\lambda$ and $\mu$ is by the ratio of
$\log(\lambda)$ to $\log(\mu)$.  The main result of this paper is that this ratio is bounded
independently of~$\phi$.

It is generally believed that most fully irreducible outer automorphisms have
different expansion factors than their inverses. In the companion paper 
\cite{HandelMosher:parageometric} we exhibit this property for a large  class of examples, the
\lq parageometric\rq\ outer automorphisms. 

Throughout this paper we say that a bound or a constant is uniform if it depends only on $n$ and not on the choice of $\phi \in \Out(F_n)$. Suppose that $A,B > 0$. If $A/ B$ is uniformly bounded above  then we write $A \prec B$ and if $A/ B$ is uniformly bounded below then we write $A \succ B$.  If $A \prec B$ and $A \succ B$ then we write $A \sim B$.

\begin{thm} \label{fully irreducible} Suppose that $\phi \in \Out(F_n)$ is fully
irreducible, that $\lambda$ is the expansion factor of $\phi$ and that $\mu$ is the expansion
factor of $\phi^{-1}$.  Then $\log(\lambda) \sim \log(\mu)$.
\end{thm}

Associated to an arbitrary outer automorphism $\phi \in \Out(F_n)$ there is a finite, indexed set of expansion factors with values in $(1,\infty)$, in terms of which we can formulate a generalization of Theorem~\ref{fully irreducible}. We recall the definition of these expansion factors from sections 3.1---3.3 of \BookOne, expressed using the language of relative train track maps and attracting laminations. 

To each $\phi \in \Out(F_n)$ there is associated (Definition 3.1.5 of \BookOne) a finite set of attracting laminations $\L(\phi)$. If $g \from G \to G$ is a relative train track representative of $\phi$, then $\L(\phi)$ contains one element $\Lambda^+$ for each subgraph $H \subset G$ which is an exponentially growing aperiodic stratum for some power of $g$; generic leaves of $\Lambda^+$ are obtained by iterating edges of $H$ and passing to a limit.  The expansion factor for the action of $\phi$ on $\Lambda^+$ (Definition~3.3.2 of \BookOne) measures the asymptotic expansion of the number of edges of $H$ in a generic leaf segment of $\Lambda^+$ under iteration.  By Proposition~3.3.3 of \BookOne, this expansion factor equals the Perron-Frobenius eigenvalue of the submatrix associated to $H$ of the transition matrix for $g$.   The set of attracting laminations $\L(\phi)$ and the expansion factor associated to each element of $\L(\phi)$ are independent of the choice of  relative train track representative. 

There is a natural bijective pairing between the expansion factors of $\phi$ and those of $\phi^\inv$, defined as follows. Associated to each $\Lambda^+ \in \L(\phi)$ is a free factor of $F_n$, called the \emph{supporting free factor} of $\Lambda^+$; this is the smallest free factor on which each leaf of $\Lambda^+$ is carried in the appropriate sense. Distinct elements of $\L(\phi)$ have distinct supporting free factors. As proved in \BookOne\ Lemma~3.2.4, there is a natural bijection between $\L(\phi)$ and $\L(\phi^\inv)$, where $\Lambda^+ \in \L(\phi)$ corresponds to $\Lambda^- \in \L(\phi^\inv)$ if and only if $\Lambda^+,\Lambda^-$ have the same supporting free factor. Thus, the expansion factors of $\phi$ and of $\phi^\inv$ can be indexed by the same set, namely the set of supporting free factors of their attracting laminations; this common indexing defines the desired bijective pairing. 

Throughout this paper, we shall denote the paired expansion factors of $\phi,\phi^\inv$ as a pair of sequences $(\lambda_i)$, $(\mu_i)$, respectively, where $i$ varies over the same finite index set. 

For completeness sake, in section~\ref{prelim} we shall give a self-contained definition of expansion factors and the bijective pairing between inverse expansion factors, formulated entirely in terms of filtered topological representatives, without recourse to the overhead needed to define relative train track maps and attracting laminations.

The natural generalization of Theorem~\ref{fully irreducible} holds in the general case.

\begin{thm} \label{reducible} Given $\phi \in \Out(F_n)$, if $(\lambda_i)$, $(\mu_i)$ are the paired expansion factors of $\phi,\phi^\inv$, respectively, then $\log(\lambda_i) \sim \log (\mu_i)$ for all $i$.  
\end{thm}

\subparagraph{Methods of proof.} Denote the rose with $n$ petals by $R_n$ and identify $\pi_1(R_n)$ with $F_n$.  A {\em marked graph} is a graph $G$, all of whose vertices have valence at least two, equipped with a homotopy equivalence $\rho_G: R_n \to G$ called the {\em marking}.   A homotopy equivalence $f :G \to G'$ between marked graphs {\em respects the markings} if $f \circ \rho_G$ is homotopic to $\rho_{G'}$.

Consider a homotopy equivalence $f \from G \to G'$ such that $f$ takes vertices to vertices, $f$ respects the markings, and for each edge $e$ of $G$, either $f(e)$ is a single vertex or $f|e$ is an immersion.  If $f(e)$ is a vertex, define $\abs{f(e)} = 0$. If $f|e$ is an immersion, define $\abs{f(e)}$ to be the edge length of  the immersed path $f(e)$;  equivalently, if $f|e$ is an immersion then  $f|e$ is one plus the number of points in the interior of $e$ that map to vertices.   Define $L(f) = \log(\sum_{e} \abs{f(e)})$ where the sum is taken over all edges $e$ of~$G$.   Thus $L(f)$ is the log of the \lq total edge length\rq\  of $f$.

Define $d(G,G') = \min_f\{L(f)\}$.   Theorem~\ref{fully irreducible} reduces (see
section~\ref{prelim}) to the following proposition.

\begin{prop} \label{quasi-metric}  $d(G,G') \sim d(G',G)$ for all marked graphs $G$ and $G'$ that have no valence~2 vertices.
\end{prop}

Recall  \cite{Vogtmann:OuterSpaceSurvey} that the spine $K_n$ of Outer Space has the structure of a locally
finite simplicial complex whose vertices are equivalence classes of marked graphs without valence~2 vertices, 
where $G$ and $\hat G$ are equivalent if there is a simplicial homeomorphism $h \from G \to \hat G$ that respects the markings.  Thus $d(G,G')$ descends to a well defined function on ordered pairs of vertices of $K_n$. It is easy to check that there exists $C > 0$ so that  $d(G,G) < C$ for all $G$ and so that $d(G,G'') \le d(G,G') + d(G',G'') +C$ for all $G, G'$ and $G''$.  Proposition~\ref{quasi-metric} can therefore be interpreted as the statement that $d(G,G')$ defines a \lq quasi-metric\rq\  on the vertex set of  $K_n$
and hence  a \lq quasi-metric\rq\ on $K_n$.

One question which arises is whether the `quasi-metric' $d$ is quasi-isometric to
the standard simplicial metric $d_s$ on $K_n$, which itself is quasi-isometric to any
word metric on $\Out(F_n)$. In other words, do there exist constants $L \ge 1$, $C \ge 0$
such that $\frac{1}{L} \, d_s(G,G') - C \le d(G,G') \le L \, d_s(G,G') + C$, for any two
marked graphs $G,G'$? This question is easily answered in the negative. Suppose that
$G_0$ is a marked rose with edges $e_1, \dots, e_n$ and that $f \from G_0\to G_0$ is
defined by $f(e_2) =  e_2e_1$ and $f(e_i) = e_i$ for $i \ne 2$.  Let $G_m$, $m > 0$, be
the marked graph obtained from $G_0$ by changing its marking so that  $f^m \from G_0 \to
G_m$ respects the markings. Then $d(G_0,G_m) \le \log(m)$, as opposed to the linear
growth of $d_s(G_0,G_m)$ proved in \cite{Alibegovic:translation}. This argument
implies, moreover, that $d$ is not quasi-isometric to any $\Out(F_n)$-equivariant proper,
geodesic metric on $K_n$, because all such metrics are quasi-isometric to
$d_s$. Whether a precise relation between $d$ and $d_s$ can be described remains unclear.

 We conclude this section with a discussion of the proof of Proposition~\ref{quasi-metric} and Theorems~\ref{fully irreducible} and \ref{reducible}.

A {\em filtration  $\F$} of a marked graph $G$ {\em with length $N$} is an increasing sequence of distinct subgraphs $\emptyset = G_0 \subset G_1 \subset \cdots \subset G_N = G$  such that no $G_i$ has  a valence one vertex.  For $1 \le i \le N$, the closure of $G_i\setminus {G_{i-1}}$ is called the {\em $i^{th}$ stratum} of the filtration and is denoted by either $H_i$ or $G(i,i)$.  For any  $f : G \to G'$ as above, denote $f(G_i)$ by $G_i'  \subset G'$.  If $G_i'$ has the same rank as $G_i$ for all $i$, we say that {\em $f$ respects the filtration}.  In this case  $\emptyset = G'_0 \subset G'_1 \subset \cdots \subset G'_N = G'$ is a filtration for $G'$ and $f$ restricts to a homotopy equivalence between $G_i$ and $G_i'$ for each $i$. 

 Every $f: G \to G'$ respects the unique   filtration $\emptyset = G_0 \subset G_1 =G$ of length one.  We may therefore assume that $f: G \to G'$ is equipped with a filtration $\F$ that it respects.  In the context of Theorem~\ref{reducible}, one expects to work with filtrations that have length greater than one.  It turns out that even in the context of Theorem~\ref{fully irreducible}, it is helpful to consider such filtrations.  In section~\ref{main prop} we state our main technical result,  Proposition~\ref{main technical}.  It is a generalization of Proposition~\ref{quasi-metric} that incorporates filtrations into its statement. The statement of Proposition~\ref{main technical}  also includes   properties that are required for our induction argument and so is   somewhat lengthy.  Using results from the theory of relative train track maps, we show in section~\ref{main prop} that Proposition~\ref{main technical}  implies Proposition~\ref{quasi-metric} and Theorems~\ref{fully irreducible} and  \ref{reducible}.  

For the sake of simplicity we focus our remaining remarks  on Proposition~\ref{quasi-metric}.

A homotopy equivalence $f \from G \to G'$ as above can be  factored   via \cite{Stallings:folding} as  
$$
G = K^0 \xrightarrow{p_1} K^1 \xrightarrow{p_2} \cdots   \xrightarrow{p_k} K^k \xrightarrow{\theta}  G'
$$
where each $p_i$ is a fold and $\theta$ is a homeomorphism.  Each  $p_k$ and $\theta$ in this {\em folding sequence} has a natural homotopy inverse.  Composing these in reverse order defines a homotopy inverse  $g \from G' \to G$ for $f : G \to G'$.  If $L(g) \prec  L(f)$ then we say  that the folding sequence has a  controlled inverse.

We  prove that every   $f :G \to G'$ as above has a folding sequence with a controlled inverse.   The proof is by downward induction on the length of a filtraton $\F$ respected by $f$.  Thus we first consider those $f$ that respect a filtration of maximal  length $2n-1$ and work our way down to those $f$ that respect only the unique filtration of length one.  

Assume that $\F$ has been specified.  A fold $p$  identifies points in a pair of edges.  If both  edges belongs to $G_j$ and if $p$ restricts to a simplicial homeomorphism of $G_{j-1}$  then we say that $p$ is supported on  $H_j$.  In Lemma~\ref{one stratum suffices}, we show that  it is sufficient to consider only those $f$  for which there exists $1 \le j \le N$ and a folding sequence  in which every fold is  supported on  $H_j$.  We say that such an $f$ is supported on  a single stratum (with respect to $\F$).

If $\F$ has maximal length and $f$ is supported on a single stratum,  then every folding sequence for $f$ has a controlled inverse.  The main steps in proving this are Lemma~\ref{short} and Lemma~\ref{j stratum}.  

At the other extreme, suppose that $\F$ has length one. There being only one stratum there is nothing gained by assuming that $f$ is supported on a single stratum.  In this case we begin with a randomly chosen folding sequence for $f$ with notation as above.  For $1 \le a \le b \le k$, define $f_{a,b} = p_b \circ \cdots \circ p_a   \from K^{a-1} \to K^b$ and say  that  $f_{a,b}$ is reducible if it  respects a filtration of length greater than $1$.   Let $a_1$ be the largest value for which $f_{1,a_1}$ is reducible.  If $a_1 < k$, consider the factorization $f = f_{a_1+2,k}p_{a_1+1}f_{1,a_1}$.   Repeating this operation with  $f_{a_1+2,k}$ replacing $f$ leads to decomposition of  $f$ as an alternating concatenation of maximal length reducible maps and single folds.     As our proof is by downward induction on the length of the filtration, we may assume that there are folding sequences with controlled inverse for each of  the maximal length reducible maps.  Together with the single folds, they define a  folding sequence for $f$ which we prove has a controlled inverse.   

For filtration lengths between the minimum and the maximum, a more general notion of reducibility is required.  This is addressed in section~\ref{irreducibility}.

\paragraph{Remark.} Although the techniques of this paper share some concepts with the techniques of relative train track maps, particularly the use of filtrations,   this paper is almost entirely independent of the theory of relative train track maps. The sole exceptions are as follows. While the statement of Lemma~\ref{train tracks} is independent of relative train track maps, the proof is an extended exercise in the methods of relative train track maps.  
Also, we will need to quote results on relative train track maps to deduce Theorem~\ref{fully irreducible} from Proposition~\ref{quasi-metric} and to justify our definitions of expansion factors in section~\ref{prelim}.

\section{Preliminaries}\label{prelim}
 All maps $h \from G \to G'$ are assumed to be homotopy equivalences of marked graphs of a fixed rank $n$. We say that $h \from G \to G'$ is {\it simplicial} if for each edge $E$ of $G$, either $h(E)$ is a vertex or $h(E)$ is a single edge  of $G'$ and $h|\Int(E) $ is injective. We assume throughout that $h \from G \to G'$  take vertices to vertices and is simplicial with respect to some subdivision of $G$. 
 
\paragraph{Filtrations.} By a {\it weak filtration $(G_i)$ of a marked graph $G$ with length $N$} we mean a properly nested sequence of subgraphs
$$\emptyset = G_0 \subset G_1 \subset \cdots \subset G_N = G
$$ 
some of which may have valence one vertices.
A \emph{filtration} is a weak filtration in which no $G_i$ has a valence one vertex, and thus all components of $G_i$ are non-contractible and $G_k$ is not homotopy equivalent to $G_{i-1}$. For sets $A \subset X$ we use the notation $X\setminus A$ for the complement of $A$ in $X$; in general $X \setminus A$ is not closed in $X$.  We denote $\closure(G_b \setminus G_{a-1})$ by $G(b,a)$.  In the special case that $b= a$ we say that  $G(b,b) = \closure(G_b \setminus G_{b-1})$ is the {\it $b$-stratum} of $G$.  A marked graph with a (weak) filtration is called a {\it (weak) filtered marked graph}. Our standing notation will be that the filtration on a graph $G$ has filtration elements $G_i$.  

Given a filtration $(G_i)$ of length $N$ of a rank $n$ marked graph $G$, the length $N$ satisfies the bound $N \le 2n-1$, because for each pair of properly nested filtration elements $G_i \subset G_{i+1}$, either $\rk(H_1(G_{i+1})) > \rk(H_1(G_i))$ or $\rk(H_0(G_{i+1})) <  \rk(H_0(G_i))$. 

Suppose that $G$ and $G'$ are marked graphs with weak filtrations of length $N$.  A homotopy equivalence $h\from G \to G'$ {\it respects the weak filtrations}  if:
\begin{itemize}
\item   $h|G_i \from G_i \to G_i'$  is a homotopy equivalence for each $1 \le i \le N$. 
\end{itemize} 
If in addition the weak filtrations are filtrations, so that no $G_i$ or $G'_i$ has any valence one vertices, then $G_i' = h(G_i)$. 

Suppose that $G$ and $G'$ are marked graphs with weak filtrations of length $N$ and that $h \from G \to G'$ respects the weak filtrations. Suppose also that $\{E_1,\ldots,E_r\}$ and  $\{E'_1,\ldots,E'_s\}$ are  the edges of $G$ and $G'$  numbered so that the  edges in $G_{i+1}$ have larger indices than edges in $G_i$ and similarly for $G'$. A path $\sigma \subset G$ with its orientation reversed is denoted  $\bar \sigma$.   Since $h$ is simplicial with respect to a subdivision of $G$, we can identify $h(E_k)$ with a (possibly trivial) word in $\{E_j', \bar E_j' \from 1 \le j \le s\}$ called {\it the edge path associated to $h(E_k)$}. We allow the possibility that $h|E_k$ is not an  immersion or equivalently that  the edge path associated to $h(E_k)$ is not a reduced word.  {\it For the remainder of the paper we identify the image of an edge with its edge path.}  The {\it transition matrix associated to $h\from G \to G'$} is the matrix $M(h)$ whose $jk^{th}$ entry is the number of times that $E_j'$ or $\bar E_j'$ occurs in  $h(E_k)$.  For $ 1 \le r \le s\le N$ let $M_{rs}(h)$ be the  submatrix of $M(h)$ corresponding to the edges of $G(s,r)$  and the edges of  $G'(s,r)$.  

It is important that we bound the size of transition matrices.  Since $n$ is fixed, we can bound the number of edges in $G$ and $G'$, and hence the size of the transition matrix, by bounding the number of valence two vertices in $G$ and $G'$.  There are several steps in the proof in which we factor a map.  The following definition and lemma are used to control the number of valence two vertices of the intermediate graphs.

The link $lk(v)$ of a  vertex $v$ in  $G$ is the set of oriented edges with $v$ as initial endpoint.  For  $h \from G \to G'$ let $lk_h(v) \subset lk(v)$ be those oriented edges whose image under $h$ is non-trivial.   There is an induced  map $Dh\from lk_h(v) \to lk(h(v))$ that sends $E$ to the first edge of $h(E)$.  We denote by $T(h,v)$  the {\it cardinality of the image $Dh(lk_h(v))$}. We record the following obvious properties of $T(h,v)$.

\begin{lemma} \label{gates} \quad
\begin{itemize}
\item  The valence of $v$ in $G$ is greater than or equal to $T(h,v)$. 
\item If  $h_1 \from G \to K$ and $h_2 \from K \to G'$ then $T(h_1,v),T(h_2,h_1(v)) \ge
T(h_2 \circ h_1,v)$.
\end{itemize}
\end{lemma}

\paragraph{Matrix lemmas.}
Our definitions of $A \prec B$ and $A \succ B$ require that $A,B > 0$.  Many of our comparisons involve $\log(C)$ for some $C \ge 1$.  To remove $C=1$ as a special case we define 
$$
\mlog(C) = \max\{1, \log(C)\}.
$$  

For any matrix $M$ denote the largest coefficient of $M$ by $LC(M)$ and   the sum of all entries in $M$ by $L(M)$.   The following lemma is obvious but useful.

\begin{lemma}  \label{direction easy} Suppose that $M_1$ and $M_2$ are non-negative integer matrices with  at most $\alpha$ rows and columns. Then $LC(M_1M_2) \le \alpha LC(M_1)LC(M_2)$ and $LC(M_1) \le L(M_1) \le \alpha^2 LC(M_1)$. 
\end{lemma} 

We will often work with non-square matrices.  For this and other reasons it is more convenient to think about largest coefficients than   about Perron-Frobenius eigenvalues.   We recall   the following well known fact. 

\begin{lemma} \label{largest} Suppose that   $M$ is a non-negative integral irreducible matrix and that $\lambda$ is its Perron-Frobenius eigenvalue.  If the size of $M$ is  uniformly bounded  then $ \mlog(\lambda) \sim  \mlog(LC(M)) \sim \mlog(L(M))$.
\end{lemma} 
\proof    Suppose that all matrices being considered have at most $\alpha$ rows and columns. 
Since $\lambda \ge 1$ is an eigenvalue with positive eigenvector for $M$, $\lambda \le \alpha \cdot LC(M)$.   On the other hand, every non-zero coordinate of $M^{\alpha}$ is greater than or equal to $LC(M)$.  This implies that $\lambda^{\alpha} \ge LC(M)$.  This proves the first comparison. The second comparison follows from Lemma~\ref{direction easy}.  
\endproof

\paragraph{Reduction of Theorem~\ref{fully irreducible} to Proposition~\ref{quasi-metric}.}
 
 A homotopy equivalence $f \from G \to G$ of a marked graph determines  an outer automorphism of $\pi_1(G)$ and so an element $\phi \in \Out(F_n)$.  We say that {\it $f \from G \to G$ represents $\phi$}.

\begin{proof}[Proposition~\ref{quasi-metric}  implies Theorem~\ref{fully irreducible}]  
By Theorem 1.7 of  \cite{BestvinaHandel:tt} there exists $f \from G \to G$ representing $\phi$  such that $M(f)$ is irreducible and has Perron-Frobenius eigenvalue~$\lambda$.   Lemma~\ref{largest} and Proposition~\ref{quasi-metric} imply that  there is a homotopy inverse $g \from G \to G$ for $f$ such that $\mlog(LC(g)) \prec \log (\lambda)$.  Remark 1.8 of \cite{BestvinaHandel:tt} implies that some submatrix of $M(g)$ is irreducible with Perron-Frobenius eigenvalue greater than or equal to $\mu$.  Another application of Lemma~\ref{largest} shows that $\log(\mu) \prec \log(\lambda)$.  Symmetry
completes the proof.  
\end{proof}

\section{The expansion factors of an outer automorphism.} 
In this section we define the expansion factors of an outer automorphism $\phi$, and the bijective pairing between the expansion factors of $\phi$ and of $\phi^\inv$. We shall state these definitions in a self-contained manner, solely in terms of filtered marked graphs and free factor systems for $F_n$, although several results from \cite{BestvinaHandel:tt} and \cite{BFH:TitsOne} regarding relative train track maps and attracting laminations will be quoted in order to prove that these definitions make sense, and that they agree with the definitions given in the introduction.

We also state and prove Lemma~\ref{train tracks}, which collects in one statement several properties of relative train track maps that will be useful in what follows.

The reader who is interested only in Theorem~\ref{fully irreducible} can skip this section entirely. The reader who is interested in Theorem~\ref{reducible} and who has a good understanding of the concepts of relative train track maps, attracting laminations, and expansion factors can just read the statement of Lemma~\ref{train tracks} and skip the rest of this section.

\paragraph{Definition of expansion factors.} We now formulate the definition of the expansion factors of an arbitrary outer automorphism $\phi \in \Out(F_n)$.
 
Let $f \from G \to G$ be a homotopy equivalence of a marked graph. Among all weak filtrations $(G_i)$ that $f$ respects there is a unique maximal one, and it has the property that the transition matrix $M_{ii}(h)$ for the $i$-stratum $G(i,i)$ is either an irreducible matrix or a zero matrix (\BH\ section~5). If $M_{ii}(h)$ is irreducible and has Perron-Frobenius eigenvalue $\lambda_i>1$ then we say that $G(i,i)$ is an \emph{exponentially growing} or \emph{EG}-stratum. Let $\Gamma(f)$ denote the sequence of Perron-Frobenius eigenvalues of EG-strata written in decreasing order $\lambda_{i_1} \ge \cdots \ge \lambda_{i_R}$, where $\{i_1,\ldots,i_R\}$ is an enumeration of the indices of the EG-strata  (see pages 32--33 of \cite{BestvinaHandel:tt} where  $\Gamma(f)$ is denoted $\Lambda(f)$). Although the sequence $\Gamma(f)$ may already have repeated entries, the multiplicity is not quite correct for defining the expansion factors of $\phi$. The correct multiplicity is determined as follows. For each irreducible square matrix $M$ there exists a unique integer $p>1$, the \emph{multiplicity} of $M$, such that some positive power $M^j$ is a block diagonal matrix with $p$-blocks each having all positive entries. Let $p_i$ be the multiplicity of $M_{ii}(f)$; if $p_i=1$ then $G(i,i)$ is called an \emph{EG-aperiodic} stratum. Let $\hat\Gamma(f)$ be the sequence obtained from $\Gamma(f)$ where $\lambda_{i}$ is replaced by $p_{i}$ equal entries of the same value as $\lambda_{i}$. To put it another way, letting $p$ be the least common multiple of the multiplicities of all EG-strata of $f$, then each EG-stratum $G(i,i)$ of $f$ of multiplicity $p_i$ breaks into $p_i$ EG-strata of $f^p$ of multiplicity $1$, and so $\hat\Gamma(f)$ is obtained from $\Gamma(f^p)$ by taking the $p^{\text{th}}$ root of each entry of $\Gamma(f^p)$. Note that the operation of replacing $\Gamma(f)$ by $\hat\Gamma(f)$ preserves lexicographic order.

Now define $\hat\Gamma(\phi)$ to be the lexicographic minimum of the set 
$$\{\hat\Gamma(f) \suchthat \text{$f \from G \to G$ represents $\phi$} \}
$$
To prove that the lexicographic minimum exists, the subset where $G$ has no valence~2 vertices is cofinal with respect to the ``greater than'' relation, according to \BH\ Lemma~5.4. But that subset is well-ordered, because if $G$ has no valence~2 vertices then it has at most $3n-3$ edges, so $f$ has at most $3n-3$ strata, and so $\Gamma(f)$ has at most $3n-3$ entries. The indexed set $\hat\Gamma(\phi)$ is defined to be the expansion factors of $\phi$.

To check that this definition of the expansion factors agrees with the one given in the introduction, we quote Theorem~5.12 from \BH\ to conclude that $\Gamma(f)$ is minimized by some relative train track representative $f_0 \from G_0 \to G_0$, and so $\hat\Gamma(f)$ is minimized as well by $f_0$. We then pass to a power of $\phi$ so that each EG-stratum of $f_0$ is aperiodic; the effect is to take the $p$th power of each element of $\hat\Gamma(f_0)$, and the $p$th power of the expansion factor of any each expanding lamination in $\L(\phi)$. Next, consider \emph{any} relative train track representative $f \from G \to G$ of $\phi$. By Lemma~3.1.14 of \BookOne, each EG-stratum of $f$ is aperiodic. Next we quote Definition 3.1.12 of \BookOne\ to conclude that the attracting laminations of $\phi$ correspond bijectively with the EG-strata of $f$. Finally we quote Proposition~3.3.3 from \BookOne\ to conclude that the expanding lamination associated to an aperiodic EG-stratum of $f$ has expansion factor equal to the Perron-Frobenius eigenvalue of that stratum.

The results just quoted also show that for \emph{every} relative train track representative $f \from G \to G$ of $\phi$ we have $\hat\Gamma(f) = \hat\Gamma(\phi)$.  

\paragraph{The bijective pairing between $\hat\Gamma(\phi)$ and $\hat\Gamma(\phi^\inv)$.} A \emph{free factor} of $F_n$ is a nontrivial subgroup $A \subgroup F_n$ such that $F_n = A * B$ for some subgroup $B$. Let $[A]$ denote the conjugacy class of $A$. Define a partial ordering on free factor conjugacy classes, where $[A] \sqsubset  [A']$ if $A$ is conjugate to a free factor of~$A'$. Note that if $F=A*A'*B$ where $A,A'$ are nontrivial then $[A] \ne [A']$. Given a free factorization $F_n = A_1 * \cdots * A_I * B$ where $A_1,\ldots,A_I$ are nontrivial, we say that $\{A_1,\ldots,A_I\}$ is an \emph{independent set of free factors}, and that $\A = \{[A_1],\ldots,[A_I]\}$ is a \emph{free factor system}. By convention the empty set $\emptyset$ is a free factor system. We extend the partial ordering $\sqsubset$ to free factor systems, where $\A \sqsubset \A'$ if for each $[A] \in \A$ there exists $[A'] \in \A'$ such that $[A] \sqsubset [A']$; the relation $\A \sqsubset \A'$ is \emph{proper} if $\A \ne \A'$. A \emph{free factor filtration $(\A^k)$ of length $K$} is a finite, properly nested sequence of free factor systems $\emptyset = \A^0 \sqsubset \A^1 \sqsubset \cdots \sqsubset \A^K = \{[F_n]\}$. The length $K$ satisfies the same bound $K \ge 2r-1$ for the same reason as the length of a marked graph filtration: for each proper relation $\A^k \sqsubset \A^{k+1}$, either the sum of the ranks of $\A^{k+1}$ is strictly greater than the sum of the ranks of $\A^k$, or the cardinality of $\A^{k+1}$ is strictly less than the cardinality of $\A^k$.


The group $\Out(F_n)$ acts on free factor conjugacy classes, systems, and filtrations. Given $\phi \in \Out(F_n)$, a $\phi$-invariant free factor filtration $(\A^k)$ is \emph{reduced} if the following holds: if a free factor system $\A'$ is invariant under the action of an iterate of $\phi$, and if $\A^k \sqsubset \A' \sqsubset \A^{k+1}$, then either $\A^k=\A'$ or $\A'=\A^{k+1}$. If this does not hold, so that $\A'$ exists with $\A^k \sqsubset \A' \sqsubset \A^{k+1}$ and $\A' \neq \A^k,\A^{k+1}$, then by passing to a power of $\phi$ which preserves $\A'$, we can insert $\A'$ and get a longer $\phi$-invariant free factor filtration.  Since the length of a free factor filtration is bounded above, it follows that some power of $\phi$ respects a reduced free factor filtration.

To each marked graph $G$ and each subgraph $G' \subset G$ there corresponds a free factor system $\F(G')$ as follows: if $G'$ is connected then $\F(G')$ is a singleton consisting of the (well-defined) conjugacy class of $\pi_1(G')$ in $\pi_1(G)$; and in general if $G'$ has components $G'=G'_1 \union \cdots \union G'_m$ then the (well-defined) conjugacy classes of $\pi_1(G'_1),\ldots,\pi_1(G'_m)$ in $\pi_1(G)$ collect together to form a free factor system $\F(G')$. To each filtered marked graph $\emptyset = G_0 \subset G_1 \subset \cdots \subset G_K=G$ there corresponds a free factor filtration $\emptyset \sqsubset \F(G_1) \sqsubset \cdots \sqsubset \F(G_K) = \{[F_n]\}$. If $\phi \in \Out(F_n)$ has a topological representative $f \from G \to G$ that respects the filtration of $G$, then $\phi$ preserves the corresponding free factor filtration. Conversely, the proof of Lemma~2.6.7 of \BookOne\ shows that for any $\phi$-invariant free factor filtration $(\A^k)$ there exists a filtered marked graph $(G_k)$ and a filtration respecting representative $f \from G \to G$ of $\phi$ such that $\A^k = \F(G_k)$. Moreover, we may choose $f$ to be a relative train track representative of $\phi$.

Consider now $\phi \in \Out(F_n)$. We shall define the bijective pairing between $\hat\Gamma(\phi)$ and $\hat\Gamma(\phi^\inv)$ after passing to a certain power of $\phi$, and one easily checks that the definition is independent of which power is chosen. Replace $\phi$ by a power so that there exists a $\phi$-invariant reduced free factor filtration $(\A^k)$. Choose a filtered marked graph $\emptyset = G_0 \subset G_1 \subset \cdots \subset G_K=G$ and a filtration respecting representative $f \from G \to G$ of $\phi$ such that $\A^k = \F(G_k)$; we may assume that $f$ is a relative train track representative of $\phi$, and so in particular $\hat\Gamma(f)$ realizes $\hat\Gamma(\phi)$. Since the free factor filtration is reduced, it follows that each stratum $G(k,k)$ is aperiodic.

The key result of \BookOne\ which establishes the pairing between $\hat\Gamma(\phi)$ and $\hat\Gamma(\phi^\inv)$ is Lemma~3.2.4. The statement of this lemma requires full knowledge of the attracting laminations of $\phi$ and of $\phi^\inv$. But the proof of Lemma~3.2.4 yields the following statement, which does not require knowledge of attracting laminations, and which can still be used to support the definition of the pairing, because it characterizes EG-strata purely in terms of free factor filtrations:
\begin{itemize}
\item If $f \from G \to G$ is a relative train track map whose strata are aperiodic and whose associated free factor filtration $(\A^k =\F(G_k))$ is reduced, then for each $k$, the stratum $G(k,k)$ is  \emph{not} EG if and only if the following holds: 
\begin{description}
\item[$(*)$] \emph{either} there exist two elements $[A],[A'] \in \A^{k-1}$ represented by independent free factors $A,A'$ so that 
$$\A^k = \left(\A^{k-1} \setminus \{[A],[A']\} \right) \union \{[A*A']\}
$$ 
\emph{or} there exist $[A] \in \A^{k-1}$ and $[A'] \in \A^k$ represented by free factors $A \subset A'$ so that $\rank(A')=\rank(A)+1$, and 
$$\A^k = \left( \A^{k-1} \setminus \{[A]\}\right) \union \{[A']\} 
$$ 
\end{description}
\end{itemize}
Using this statement we can define the pairing between $\hat\Gamma(\phi)$ and $\hat\Gamma(\phi^\inv)$ as follows. We are assuming that $(\A^k)$ is a reduced free factor filtration invariant under $\phi$, and so also invariant under $\phi^\inv$. Choose relative train track representatives $f \from G \to G$ for $\phi$ and $f' \from G' \to G'$ for  $\phi^{-1}$ so that $\F(G_k)=\F(G'_k) = \A^k$, for all $k=0,\ldots,K$. Applying the above proposition, it follows that for each $k=0,\ldots,K$, $G(k,k)$ is an EG-stratum of $f$ if and only if $G'(k,k)$ is an EG-stratum of $f'$, because both of these are equivalent to the failure of property $(*)$. We therefore define the pairing by requiring that the Perron-Frobenius eigenvalue of $f$ on $G(k,k)$ be paired with the Perron-Frobenius eigenvalue of $f^\inv$ on $G'(k,k)$.

The fact that this is well-defined independent of the choice of relative train track representatives $f$ and $f'$, and that this definition agrees with the one given in the introduction, is again a direct consequence of the proof of Lemma~3.2.4, which shows that the attracting lamination of $\phi$ associated to $G(k,k)$ is paired with the attracting lamination of $\phi'$ associated to $G'(k,k)$, and so the expansion factor $\lambda_k$ of $\phi$ associated to $G(k,k)$ is paired with the expansion factor $\mu_k$ of $\phi^{-1}$ associated to $G'(k,k)$.

\paragraph{A relative train track lemma.} The proof of Theorem~\ref{fully irreducible} followed from Proposition~\ref{quasi-metric} plus a modest amount of train track theory quoted from \cite{BestvinaHandel:tt}, as we saw at the end of section~\ref{prelim}. For Theorem~\ref{reducible} we will need more. The proof of the following lemma is essentially an extended exercise in relative train track theory.

\begin{lemma} \label{train tracks} Suppose that $\phi \in \Out(F_n)$ and that $\lambda$ and $\mu$ are paired expansion factors of $\phi$. Then there is a filtered marked graph $G$ with filtration $(G_k)$, a homotopy equivalence $f \from G \to G$ that respects the filtration and that represents $\phi^p$ for some $p \ge 1$,  and there is a stratum $G(j,j)$ such that:
\begin{description}
\item [(1)] $f$ restricts to an immersion on each edge of $G$.
\item [(2)] $T(f | G_i,v) \ge 2$ for all vertices $v \in G_i$ and for all $i$.
\item [(3)] $\mlog(LC(M_{jj}(f))) \sim p \log(\lambda)$.
\item [(4)] if $g \from G \to G$ is a homotopy inverse for $f \from G \to G$ that respects the filtrations then  $\mlog(LC(M_{jj}(g))) \succ p \log(\mu)$.
\end{description}
Moreover, $G$ has at most $V_n'$ vertices for some uniform constant $V_n'$. 
\end{lemma}

\begin{com} The uniformity of $V_n'$ expresses the fact that relative train track graphs can be chosen with a uniformly bounded number of valence two vertices.  One can choose $G$ so that no embedded arc in $G$ contains three consecutive valence two vertices, thereby giving an explicit value for $V_n'$ of $6n-6$.   
\end{com}

\begin{notn} Choose a uniform constant $V_n \ge 2V_n'$ such that a marked graph of rank $n$ has at most $V_n/2$ vertices of valence at least three.  Thus a marked graph of rank $n$ with at most $V_n/2$ vertices of valence two has at most $V_n$ vertices.  Choose a uniform constant $Edge_n$ such that a graph of rank $n$ with at most $V_n$ vertices has at most $Edge_n$ edges.  
\end{notn}

\paragraph{The proof of Lemma \ref{train tracks}.}  Fix $\phi \in \Out(F_n)$ and paired expansion factors $\lambda$ of $\phi$ and $\mu$ of $\phi^\inv$.

A relative train track map $f\from G \to G$ is a homotopy equivalence of a marked graph $G$ which has useful properties for computation --- explicit details can be found on page 38 of \cite{BestvinaHandel:tt} and in Theorem 5.1.5 of \cite{BFH:TitsOne}. We will not use the most technical properties listed in this theorem and will be explicit about what is necessary. 

Notation: the maximal weak filtration preserved by a relative train track map $f \from G \to G$ will be denoted $(\hat G_i)$, with strata denoted $\hat H_i$, as is conventional when considering relative train track maps, rather than $\hat G(i,i)$. Also let $\hat M_{ii}$ denote the transition matrix of the restriction of $f$ to the stratum $\hat H_i$. The ``hat'' notation $\hat G_i$ is used to distinguish this weak filtration from the (non-weak) filtration $(G_k)$ with strata $H_k = G(k,k)$ whose description is the goal of the proof of Lemma~\ref{train tracks}. Also, we must describe a particular index $j$ of an EG-stratum $H_j$ associated to $\lambda$.

Choose $p > 0$ so that Theorem 5.1.5 of \cite{BFH:TitsOne} applies to both $\phi^p$ and $\phi^{-p}$ and let $f \from G \to G$   be  a relative train track map $f \from G \to G$   representing $\phi^p$ and satisfying the conclusions of that theorem.  Item (1) of Lemma~\ref{train tracks} holds for every relative train track map (and can be arranged simply by tightening the image of edges). In particular, $Df$ is defined on every oriented edge in the link of a vertex.  
  
The strata $\hat H_i$ of the weak filtration $(\hat G_i)$ are divided into four mutually exclusive types.  
\begin{description}
\item[Exponentially growing:] $\hat H_i$ is an EG-stratum. As we have seen in section~\ref{prelim}, there is a bijection between the set of exponentially growing strata $\hat H_i$ and the set of paired expansion factors $\lambda_i,\mu_i$ of $\phi,\phi^\inv$, respectively. If $\hat H_i$ corresponds to the pair $\lambda_i,\mu_i$ then $PF(\hat M_{ii}) = \lambda_i^p$.
\item[Zero stratum:] If $\hat G_i$ has contractible components, then $\hat H_i$ is the union of these components and is called a {\it zero stratum}. In this case $f(\hat H_i) \subset \hat G_{i-1}$.
\end{description}
All the remaining strata $\hat H_i$ are single oriented edges $\hat E_i$. 
\begin{description}
\item[Fixed edge:] If $f(\hat E_i) = \hat E_i$ then both endpoints of $\hat E_i$ have valence at least two in $\hat G_{i-1}$ and  $\hat H_i$ is a {\it fixed edge}.  In some contexts fixed non-loop edges are collapsed but the resulting filtration might not be reduced so we do not do that here.   
\item[Non-exponentially growing:] Otherwise (see (ne-i) and (ne-ii) in the statement of  Theorem 5.1.5 of \cite{BFH:TitsOne})    $f(\hat E_i) =  \hat E_i  \hat u_i$ for some non-trivial closed path $\hat u_i \subset  \hat G_{i-1}$ that is immersed both as a path and as a loop. In this case the terminal endpoint of $\hat E_i$  has valence at least two in $\hat G_{i-1}$.   Strata of this fourth kind are said to be {\em \neg}.
\end{description}

We define the filtration $(G_k)$ from the bottom up with each stratum of $(G_k)$ being a union of strata of $(\hat G_i)$. The stratum $H_j$ associated to $\lambda$ will be the one that contains the weak EG-stratum $\hat H_i$ associated to $\lambda$. Since $\hat G_0 = \emptyset$,  $\hat H_1$ is either a fixed loop or is \eg\ and we define $G_1 =\hat G_1$. In the former case, (2) is clear and $j \ne 1$.  Suppose then that $\hat H_1$ is \eg\ and that $v \in \hat H_1$. Inductively define $Df^m(e) = Df^{m-1}(Df(e))$.  We say that a pair of edges  $(e,e')$ in the link of $v$ determine a {\it legal turn} if $Df^m(e)$ and $Df^m(e')$ are distinct for all $m > 0$. Clearly (2) is satisfied at $v$ if there is a legal turn at $v$.  Since $\hat H_1$ is \eg, there exists $m > 0$ and a point $x$ in the interior of an edge such that $f^m(x) = v$. The \lq link\rq\ at $x$ contains two directions and condition (RTT-3) on page 38 of \cite{BestvinaHandel:tt} implies that their images in $lk(v)$ define a legal turn. This verifies (2) for $G_1$.  If $\lambda$ corresponds to $\hat H_1$, then Condition (3) for $G_1$ follows from Lemma~\ref{largest} and the fact that $PF(\hat M_{11}) = \lambda^p$. This completes our analysis of $G_1$.

Let $l  > 1$ be the first parameter value for which the non-contractible components of $\hat G_l$ do not deformation retract to $\hat G_1$. Since a zero stratum is a union of contractible components, $\hat H_l$ is not a zero stratum. We will refer to the strata above $G_1$ and below $\hat H_l$ as {\it intermediate strata}.  Since every bi-infinite path in $G_{l-1}$ is contained $G_1$, no intermediate stratum is \eg.  Property (z-i) on page 562 of \cite{BFH:TitsOne} implies that at most one intermediate stratum is a zero stratum and if there is one, then it is $\hat H_{l-1}$ and $\hat H_l$ is \eg. 

 If $\hat H_l$ is a fixed edge then push all intermediate edges of $\hat G_l \setminus  G_1$  up the weak filtration by reordering the strata.  Thus $l = 2$ and we define $G_2 := \hat G_2$.  As in the $G_1$ case, (2) is satisfied and $j \ne 2$. Suppose next that $\hat H_l$, and hence each intermediate stratum, is \neg.  Each of these strata is a single edge  with terminal endpoint in $ G_1$.  The initial endpoint of $E_l$ is either in $ G_1$ or is shared with one other of the intermediate edges.  Push all other intermediate edges of $\hat G_l \setminus  G_1$  up the weak filtration by reordering the strata.   Thus  $l$ is two or three and $G_2 := \hat G_l$ is obtained from $G_1$ by topologically adding an arc which may be subdivided at one point.  As in the previous case (2) is satisfied and $j \ne 2$.

It remains to consider the case that  $\hat H_l$ is  \eg.  If $\hat E_i$ is the edge of an intermediate \neg\ stratum then its initial endpoint $w_i$  is fixed and so is not a vertex in a zero stratum.  If $w_i$ is not a vertex in $\hat H_l$ then we push $\hat E_i$ up the weak filtration as before.  Once this is done define $G_2 = \hat G_l$, so
each remaining intermediate stratum is either a zero stratum comprising $\hat H_{l-1}$ or an NEG stratum with initial point on $\hat H_l$ and terminal point on $G_1$.

Condition (RTT-1) on page 38 of \cite{BestvinaHandel:tt} implies that if $e, e' \in lk(v)$ where  $e$ is an edge in $\hat H_l$ and $e'$ is an edge in $\hat G_{l-1}$,  then $(e,e')$ is legal.  If the link of $v$ in $G_2$ is entirely contained in $\hat H_l$ then  we produce legal turns at $v$ as we did for the $G_1$ case.  If $v$ is in a zero stratum then (2) follows from conditions z-(ii) and (z-iii) on page 562 of \cite{BFH:TitsOne}.  This verifies (2).  It also shows that $G_2$ is $f$-invariant since no valence one vertex $w_i$ in the original $\hat G_l$ could be the image of either a vertex in $G_2$ or the interior of an edge in $G_2$.  

     We now consider (3) for $G_2$.  It is clear that the only columns of $M_{22}(f)$ that have entries greater than one are those that correspond to the edges of $\hat H_l$.   These columns record the number of times that the $f$-image $\alpha$ of an edge in $\hat H_l$ crosses an edge of $H_2$.  Denote the union of the intermediate strata by $X$ and write $\alpha$ as an alternating concatenation of subpaths $\alpha = \sigma_1 \tau_1 \sigma_2\dots$ where $\tau_i$ is contained in $X$ and $\sigma_i$ is disjoint from $X$.  Since $X$ is a forest, the number of edges in $\tau_i$ is uniformly bounded and either $\sigma_{i-1}$ or $\sigma_{i+1}$ is contained in $\hat H_l$.  It follows that the largest coefficient of $M_{22}(f)$ determined by $\alpha$ is comparable to the largest coefficient of $\hat M_{ll}(f)$ determined by $\alpha$.  Condition (3) now follows as in the previous cases.

After finitely many such steps we have defined $(G_k)$ satisfying (1) - (3). We say that a vertex of $G$ is mixed if its link intersects more than one stratum of $(G_k)$.   By the fourth item in the statement of Theorem 5.1.5 of \cite{BFH:TitsOne}, $f(v) \in \Fix(f)$ for each vertex $v$.   Redefine the simplicial structure on $G$ by eliminating from the zero skeleton any valence two vertex of $G$ that is not mixed, is not the image of a mixed vertex, and is not the image of a vertex of valence greater than two.  In this new structure, $G$ has a uniformly bounded number of vertices and so a uniformly bounded number of edges. Properties (1) - (3) are not effected by this reverse subdivision so it  remains to verify (4).

There is no loss in replacing $G$ by the component of $G_j$ that contains $H_j$. We may therefore assume that $j = N$.  We first verify (4) in the special case that $g | G_{N-1} \from G_{N-1}  \to G_{N-1}$ is a relative train track map.

By Theorem 5.1.5 of \cite{BFH:TitsOne} there is a relative train track map $g' \from G' \to G'$ representing $\phi^{-p}$, a stratum $G'_m$ of the associated filtration and a bijection between the  components $\{C_i'\}$ of $G'_m$ and the components $\{C_i\}$ of $G_{N-1}$ such that $\F(C'_i) = \F(C_i)$. Let $A_i$ be a free factor representing $\F(C'_i)=\F(C_i)$, so $[\phi^p(A_i)]=[A_i]$, and by conjugating $A_i$ back to $\phi^p(A_i)$ we obtain a well-defined outer automorphism of $A_i$ that we denote $\phi^p \restrict [A_i]$. In our special case, both $g \restrict C_i$ and $g' \restrict C_i'$ are relative train track maps for the outer automorphism $\psi_i := \phi^{-p} |  [A_i]$.  Thus $\Gamma(g\restrict C_i) = \Gamma(g' \restrict C_i')$.   Lemma 3.2.4 of \cite{BFH:TitsOne} implies that if $\Lambda_j^{\pm}$ is a lamination pair for $\phi^p$ and if $\Lambda^+_j$ is carried by $A_i$ then $\Lambda^-_j$ is carried by $A_i$. The expansion factor for the action of $\psi_i$  on any such $\Lambda^-_j$ is the same as the  expansion factor for the action of $\phi^{-p}$ on $\Lambda^-_j$.  Thus the elements of $\Gamma(g \restrict C_i)$ are the expansion factors for the action of $\phi^{-p}$ on those laminations $\Lambda^-_j$ whose paired $\Lambda^+_j$ is carried by $A_i$. It follows that the elements of $\Gamma(\phi^{-p})$ and the elements of $\Gamma(g)$ differ only in that $\mu^p$ is replaced by the finite collection of eigenvalues corresponding to irreducible weak strata in $H_N$.  At least one of these eigenvalues must be greater than or equal to $\mu^p$.  Lemma~\ref{largest} completes the proof in the special case.

We reduce to the special case as follows.  For each component $C_i$ of $G_{N-1}$, $g \restrict C_i$ determines an outer automorphism of some lower rank free group; let $h_i \from X^i \to X^i$ be a relative train track map representing this outer automorphism.  Define $X = \cup X_i$ and $h \from X \to X$ by $h|X_i =  h_i$.   There is a sequence of Whitehead moves that deform $G_{N-1}$ to $X$.  These  extend to Whitehead moves on $G$ in the obvious way to produce a deformation of $G$ to a marked graph $G^*$ with $X$ as a subgraph and with a natural bijection between the edges of $G^* \setminus X$ and the edges of $H_N$.  Let $\alpha \from G \to G^*$ be the induced marking-preserving homotopy equivalence and let $\beta\from G^* \to G$ be its homotopy inverse.   Let  $\hat g= \alpha \circ g \circ \beta \from G^* \to G^*$. Then 
\begin{itemize}
\item  $\hat g$ represents $\phi^{-p}$.
\item  $M_{NN}(\hat g) = M_{NN}(g)$. 
\item   $\hat g \restrict X$ is homotopic  to a relative train track map.
\end{itemize}

The homotopy  of the third item can be performed while maintaining the second and this
completes the proof.  \qed

\section{The Main Proposition} \label{main prop}

We will prove Proposition~\ref{quasi-metric} by working in a more general context, one that is suited to our induction technique.  We have already started laying the groundwork for this by introducing filtrations.  Our main technical result, Proposition~\ref{main technical}, is stated below.   In this section we  show that it implies Proposition~\ref{quasi-metric} and that it and Lemma~\ref{train tracks} imply Theorem~\ref{reducible}. The proof of Proposition~\ref{main technical} occupies the rest of the paper save the last section, where Lemma~\ref{train tracks} is proved.

Suppose that $G$ and $G'$ are  marked graphs  with filtrations  of length $N$ and  that  $f \from G \to G'$ respects  the filtrations.  We say that $f \from G \to G'$ is {\it supported on $G(b,a)$} if:
\begin{description}
\item [(s1)] $f|G_{a-1} \from G_{a-1} \to G'_{a-1}$ is a simplicial homeomorphism.
\item [(s2)] $f|\closure(G \setminus G_b) \from \closure(G \setminus G_b)  \to \closure(G' \setminus G'_b)$ fails to be a simplicial homeomorphism only in that it may identify distinct vertices in  $G_b$.  Moreover, $f$ maps vertices in $G \setminus G_b$ to vertices in $G' \setminus G'_b$. 
\end{description}  
In the special case that  $b = a$   we say that {\it $f$ is supported on a single stratum}.

The {\it frontier of $G_b$} is the set of vertices of $G_b$ whose link contains an edge in $G(N,b+1)$ or equivalently in some stratum higher than $b$. For $a \le b$, define {\it $\Fr(G,b,a)$} to be the set of  vertices in the frontier of $G_b$ that are not contained in $G_{a-1}$.  If $f$ is supported on $G(b,a)$ then  $\Fr(G',b,a) \subset f(\Fr(G,b,a))$ and hence     $|\Fr(G',b,a)| \le |\Fr(G,b,a)|$.   There are two ways that this inequality can be strict.  One is that $f$ identifies a pair of vertices in $\Fr(G,b,a)$ and the other is that $f$ identifies a vertex in $\Fr(G,b,a)$ with a vertex in $G_{a-1}$.

\begin{proposition} \label{main technical} Suppose that $G$ and $G'$ are  marked graphs  of rank $n$ with filtrations  of length $N$ and  that  $f \from G \to G'$ respects  the filtrations and is supported on $G(b,a)$.  Suppose further that :
\begin{itemize}
\item  $f$ restricts to an immersion on each edge.
\item  $T(f|G_i,v) \ge 2$ for all vertices $v \in G_i$ and all $i$.
\item  $G$ has at most $V_n/2$ vertices with $T(f,v) = 2$.
\end{itemize}
Then there is a homotopy inverse $g\from G' \to G$ that respects the filtrations such that \begin{description}
\item [(1)] $\mlog(LC(M_{rs}(g))) \prec \mlog(LC(M_{rs}(f)))
$ for all $1 \le r \le s \le N$. 
\item [(2)] $g|G_{a-1}'\from G_{a-1}' \to G_{a-1}$ is a simplicial homeomorphism. 
\item [(3)] if $|\Fr(G,b,a)| = |\Fr(G',b,a)|$ then $g$ is supported on $G'(b,a)$.
\end{description}
\end{proposition}

\vspace{.1in}

\noindent{\bf  Proposition~\ref{main technical} implies Proposition~\ref{quasi-metric}}  Suppose that $d(G,G') =  L(f)$ for $f \from G \to G'$.  We want to apply Proposition~\ref{main technical} with $N= a =b = r = s =1$. The first and third item in the hypothesis are satisfied so consider the second.  If $T(f,v) = 1$ then there is a vertex $v$ and an edge $e'$ in $G'$ so that the edge path $f(e)$ begins with $e'$ for every edge $e$ in $G$ incident to $v$.   There is an  obvious homotopy of $f$ that slides $f(v)$ across $e'$.  Denote the resulting map by $f' \from G \to G'$. Then $L(f') < L(f)$ in contradiction to our choice of $f$.  This verifies the second item.  Let $g\from G' \to G$ be the homotopy inverse produced by Proposition~\ref{main technical} and let $g' \from G' \to G$ be obtained from $g$ by a homotopy relative to the vertices of $G'$ that tightens the $g$-image of each edge into an immersed path.   The entries in $M(g')$ are all less than or equal to the corresponding entries of $M(g)$.  Conclusion (1) of Proposition~\ref{main technical} therefore implies that $d(G',G) \prec d(G,G')$.  Symmetry completes the proof. \qed   

\vspace{.1in}

\noindent{\bf Lemma~\ref{train tracks} and  Proposition~\ref{main technical} imply Theorem~\ref{reducible}}  Given $\phi \in \Out(F_n)$ and  a lamination pair $\Lambda^{\pm}$  for $\phi \in \Out(F_n)$ with expansion factors $\lambda$ and $\mu$, let $f \from G \to G$, $p$ and $G(j,j)$ be as in Lemma~\ref{train tracks}.  Apply Proposition~\ref{main technical} to produce a homotopy inverse $g \from G \to G$ such that 
$$
\mlog(LC(M_{jj}(g))) \prec \mlog(LC(M_{jj}(f))).
$$  
By Lemma~\ref{train tracks}  
$$
p\log \mu \prec \mlog(LC(M_{jj}(g)))
$$
and  
$$
\mlog(LC(M_{jj}(f))) \sim p \log \lambda.
$$
Thus $ \log \mu \prec \log \lambda$.  The symmetric argument shows that $\log \lambda \prec \log \mu$ and hence that    $\log \mu \sim \log \lambda$ as desired. \qed

\section{Folding}

We ultimately construct  $g \from G' \to G$   by factoring $f \from G \to G'$ as a sequence of elementary maps (folds and homeomorphisms) and inverting each one.  Most of this is done recursively.  In this section we  construct $g$ explicitly in two basic cases.

We begin with the easiest case. 

\begin{lemma}\label{homeo}  Proposition~\ref{main technical} holds for a 
 homeomorphism $f \from G \to G'$.  
\end{lemma}

\proof     If $f$ is simplicial then $g = f^{-1}$.  In general, $f$ is a  finite subdivision of $G(b,a)$
followed by a simplicial homeomorphism.  The inverse $h$ of a finite subdivision is a simplicial map that collapses finitely many edges  and we define $g = h f^{-1}$. Thus $g$ is supported on $G'(b,a)$ and $LC(M(g)) = 1$. 
\endproof

We set notation for the folding method of Stallings \cite{Stallings:folding} as follows.  

\begin{notn}   \label{foldingNotation} Suppose that $f\from G \to G'$ is a homotopy equivalence of marked graphs and that $f$ restricts to an immersion on each edge.  Suppose further that 
\begin{description}
\item  [(F1)] $\hat e_1$ is a non-trivial initial segment of  an edge $e_1$ in $G$ with initial vertex $v_0$ and terminal vertex $v_1$.
\item  [(F2)]  $\hat e_2$ is a non-trivial  initial segment of  an edge $e_2$ in $G$ with initial vertex $v_0$ and terminal vertex $v_2$.  
\item  [(F3)]   $f(\hat e_1) = f(\hat e_2)$ is the {\it maximal} common initial subpath of $f(e_1)$ and $f(e_2)$. 
\end{description}

Then $f = f_1 \circ p$ where $p \from G \to G^*$ is the quotient map that identifies $\hat e_1$ with $\hat e_2$  and $f_1 \from
G^* \to G'$ is the induced homotopy equivalence.  The map $p$ is called  {\it a fold}.  

{\em Assume now that $f$ is supported on $G(b,a)$}.  After interchanging $e_1$ and $e_2$ if necessary we may assume that
\begin{description}
\item [(F4)] $e_1 \subset G(b,a)$ and $e_2 \subset G_b$.
\item  [(F5)] If $e_2$ is  contained in $G_{a-1}$ then $\hat e_2 = e_2$.
\item  [(F6)] If $\hat e_1 = e_1$ then $\hat e_2 = e_2$ and    $v_1 \not \in G_{a-1}$; moreover if $v_1 \in \Fr(G,b,a)$ then either $v_2 \in \Fr(G,b,a)$ or $v_2 \in G_{a-1}$. 
\end{description}
\end{notn}

There are several places in the proof where we factor $f \from G \to G'$ as a finite composition.    The following lemma is used to equip the intermediate graphs with filtrations.

\begin{lemma} \label{push forward} Suppose that $G$ and $G'$ are  filtered marked graphs, that  $f \from G \to G'$ respects  the filtrations and that $G \xrightarrow{p} G^* \xrightarrow{f_1} G'$ are as
in Notation~\ref{foldingNotation}. Suppose further that  
\begin{enumerate}
\item  $f$ restricts to an immersion on each edge.
\item  $T(f|G_i,v) \ge 2$ for all vertices $v \in G_i$ and for all $i$.
\item  $G$ has at most $V_n/2$ vertices with $T(f,v) = 2$.
\end{enumerate}
Then     $G_i^* = p(G_i)$ defines a filtration ${\cal F}^*$ on $G^*$ such that    $p$ and $f_1$ respect the filtrations.   Moreover, $p$ and $f_1$ satisfy (1) - (3).
\end{lemma}

\proof  This is immediate from the definitions except perhaps for showing that $G^*$ has at most $V_n/2$ vertices with $T(f_1,v) < 3$.  If every vertex of $G^*$ is the image of a vertex of $G$ then this follows from Lemma~\ref{gates}.  If there is a new vertex $w^*$ of $G^*$ then   $\hat e_1 \ne e_1$ and  $\hat 
e_2 \ne e_2$.  The vertex $w^*$ is the image of the terminal endpoint of  $\hat e_1$ and also the image of the terminal endpoint of  $\hat e_2$. There are three edges of $G^*$ incident to $w^*$. We label them $\bar e^*, e_1^*$ and $e_2^*$ where $p(e_i) = e^* e_i^*$ for $i=1,2$.  The maximality condition on $\hat e_1$ and $\hat e_2$ guarantees that $Df_1(e_1^*),Df_1(e_2^*)$ and $Df_1(e^*)$ are all distinct. Thus $T(f_1,w^*) = 3$.   
\endproof

The filtration ${\cal F}^*$ on $G^*$ defined in Lemma~\ref{push forward} is called  {\it the pushed forward filtration}. 

\vspace{.1in}

We now explicitly invert a single fold $p \from G \to G^*$ that is supported on $G(b,a)$.  The construction depends not only on $p$ but also on $a$ and $b$.

\begin{defn} \label{explicit inverse}  There are three cases to consider.

\noindent{\bf (Case 1:  $\hat e_1 \ne  e_1$ and $\hat e_2 = e_2$)}   In this case $G^*$ is obtained from $G$ by replacing  $e_1$ with an edge $e_1^*$ that has initial endpoint $v_2$ and terminal endpoint $v_1$.   Identify   $G \setminus \{e_1\}$ with    $G^* \setminus \{e_1^*\}$.  Then $p(e_1) = e_2 e_1^*$ and $p$ is the identity on all other edges.  
The homotopy inverse $q \from G^* \to G$  is the identity on edges other than $e_1^*$ and satisfies  $q(e_1^*) = \bar e_2 e_1$.  We may view $q$ as the fold of an initial segment of $e_1^*$ with all of $\bar e_2$. In this case  $q$ is supported on  $G^*(b,a)$.    

\noindent{\bf (Case 2 :  $\hat e_1 \ne  e_1$ and  $\hat e_2 \ne  e_2$)}    By (F5),   $e_1,e_2 \subset G(b,a)$.    In this case $G^*$ is obtained from $G$ by \lq blowing up\rq\  $v_0$.  More precisely, we add a vertex $w^*$ and an edge $e^*$ connecting $v_0$ to $w^*$; the edge $e_i$, $i = 1,2$  is replaced by an edge $e_i^*$ that initiates at $w^*$ and terminates at $v_i$. Identify   $G \setminus \{e_1, e_2\}$ with    $G^* \setminus \{e^*,e_1^*, e_2^*\}$.   Then $p(e_i) = e^* e_i^*$ for $i =1,2$ and $p$ is the identity on all other edges.   The link  $lk(w^*) = \{\bar e^*,  e^*_1, e^*_2\}$ is contained entirely in $G^*(b,a)$.
The homotopy inverse $q \from G^* \to G$  collapses $e^*$ to $v$, maps $e_i^*$ to $e_i$ for $i =1,2$ and is the identity on all other edges.  As in the previous case,  $q$ is supported on $G^*(b,a)$.

\noindent{\bf (Case 3:   $\hat e_1 =  e_1$)}.  In this case (F6) applies.  The graph $G^*$ is obtained from $G \setminus \{e_1\}$ by identifying $v_1$ and $v_2$; we label the resulting vertex $v^*$.  The quotient map $p$ satisfies $p(e_1) = e_2 $ and  is the identity on all other edges except that it identifies the vertices $v_1$ and $v_2$. 

By (F6), $|\Fr(G^*,b,a)| < |\Fr(G, b,a)|$ if and only if $v_1 \in \Fr(G,b,a)$. 
The map $p$ induces an injection (also called $p$) of $lk(v_1) $ into $lk(v^*)$.  If $e_k^* = p(e_k)$ for $e_k \in lk(v_1) \setminus \bar e_1$    then $q(e^*_k) = \bar e_2 e_1 e_k$.  All other edges are mapped by the identity.  Thus $q$  is supported on  $G^*(b,a)$ if and only if  $v_1 \not \in \Fr(G,b,a)$.
\end{defn} 

\begin{lemma} \label{one fold}  Proposition~\ref{main technical} holds for single folds $p \from G \to G^*$.  
\end{lemma}
\proof  We have constructed a homotopy inverse $q \from G^* \to G$ that satisfies (2) and (3) and has $LC(M(q))=1$.
\endproof

\section{Bounded Products}  

The following lemma states that the set of maps satisfying Proposition~\ref{main technical} is closed under composition of a uniformly bounded number of factors.    

\begin{lemma} \label{product}  To prove Proposition~\ref{main technical} it suffices to prove that there exists $k > 0$ such that every $f$ satisfying the hypotheses of Proposition~\ref{main technical} factors as   $G = K^0 \xrightarrow{f_1} K^1 \xrightarrow{f_2} \cdots
\xrightarrow{f_k} K^k = G'$  where the filtration  on $K^j$ is defined, inductively, as the pushforward by 
$f_j$ of the filtration on $K^{j-1}$ and where 
each $f_j \from K^{j-1} \to K^j$ satisfies the conclusions of Proposition~\ref{main
technical}.
\end{lemma} 

\proof   Let $g_j \from K^j \to K^{j-1}$ be a homotopy inverse for $f_j \from K^{j-1} \to
K^j$  that satisfies (1) - (3) of Proposition~\ref{main technical} and define $g = g_1
\circ \dots \circ g_k$. Properties (2) and (3) are preserved under composition so are
satisfied by $g$.

By Lemma~\ref{direction easy}  
$$
LC(M_{rs}(g)) \le (Edge_n)^{k-1}\prod_{j=1}^k  LC(M_{rs}(g_j))
$$
and so by (1) of Proposition~\ref{main technical} 
$$
\mlog(LC(M_{rs}(g))) \prec  \sum_{j=1}^k \mlog(LC(M_{rs}(g_j))) \prec \sum_{j=1}^k \mlog(LC(M_{rs}(f_j))). 
$$
  For each $j=1,\ldots,k$ 
there   is an edge $e_{j-1} \subset K^{j-1}(s,r)$ and an {\it oriented} edge $e_j^* \subset K^{j}(s,r)$ that occurs at least  $LC(M_{rs}(f_{j}))/2$ times in $f_{j}(e_{j-1})$.  We may divide $e_{j-1}$ into $(LC(M_{rs}(f_j))/2)-1$ subintervals each of which has image an immersed loop that contains $e_{j}^*$.  The image of each of these immersed loops under $f_k \circ \cdots \circ f_{j+1}$   contains an edge of $G(s,r)$.  This proves that 
$$
\mlog(LC(M_{rs}(f_j))) \prec \mlog(LC(M_{rs}(f)))   
$$
and hence that
 $$
\mlog(LC(M_{rs}(g))) \prec k \mlog(LC(M_{rs}(f))) \prec \mlog(LC(M_{rs}(f)))
$$
as desired.
\endproof   

  The first application of Lemma~\ref{product} is to prove a special case
 of Proposition~\ref{main technical}.  In the context of
  Notation~\ref{foldingNotation} F4--F6, where the fold $p \from G \to G^*$ is supported on
  $G(b,a)$, if $e_2 \subset G_{a-1}$, then we say that $p \from G \to G^*$
  {\it folds into lower strata}.

\begin{lemma} \label{pg}  Proposition~\ref{main technical} holds under the additional hypothesis that $f$ factors as $G = K^0 \xrightarrow{p_1} K^1 \xrightarrow{p_2} \cdots
\xrightarrow{p_l} K^l = G'$ where each $p_i$  folds into  lower strata (with respect to the pushed forward filtrations).
\end{lemma}
\proof    The folds $p_i$ are necessarily in the first or third cases of Definition~\ref{explicit inverse}.   The number of edges in $K^i(b,a)$  is less than or equal to the number of edges  in $K^{i-1}(b,a)$  with equality if and only if $p_i$ is a case one fold.  It follows that there is a uniform bound to the number of case three folds.  By Lemma~\ref{product} and Lemma~\ref{one fold} we are reduced to the case that each $p_i$ is a case one fold.  

    There is a natural bijection between the edges of $K^{i-1}$ and the edges of $K^i$.  With respect to this bijection, $M(p_i)$ has ones on the diagonal and has exactly one non-zero off diagonal entry in a column corresponding to an edge in $K^{i-1}(b,a)$ and a row corresponding to an edge of $K^i_{a-1}$.  Recall that in the indexing convention for transition matrices, the indices for edges in $K^{i-1}(b,a)$ are all larger than the indices for edges in $K^i_{a-1}$.  The $M(p_i)$'s therefore commute.  

     Let $q_i \from K^i \to K^{i-1}$ be the inverse of $p_i$ as specified in Definition~\ref{explicit inverse} and note that $M(p_i) = M(q_i)$. Define $g = q_l \circ \ldots \circ q_1$.  Then  $M(f) = M(p_l)\cdot M(p_{l-1})\cdots M(p_1) =    M(q_1)\cdot M(q_2)\cdots M(q_l) = M(g)$.  This implies (1) of Proposition~\ref{main technical}. 
Conditions (2) and (3)  follow from Lemma~\ref{one fold} and the observation that these conditions are preserved under composition.
   
\endproof

The second application of Lemma~\ref{product} is a reduction to the case that $f$ is supported on $G(j,j)$ for some $1 \le j \le N$.

\begin{lemma} \label{one stratum suffices} It suffices to prove Proposition~\ref{main technical} with the additional hypothesis that  $f$ is supported on a single stratum.
\end{lemma}

\proof    We claim that $f$ factors as 
$$
G = K^{a-1} \xrightarrow{f_a} K^a \xrightarrow{\theta_a} K^{a+1} \xrightarrow{f_{a+1}} K^{a+2} \xrightarrow {\theta_{a+1}}  
\cdots \xrightarrow{f_{b}} K^{2b-a} \xrightarrow{\theta_{b}} K^{2b-a+1} = G'
$$
  where for   $a \le j \le b$:
\begin{itemize}
\item  $f_j \from K^{2j-a-1} \to K^{2j-a}$ is a  homotopy equivalence that respects the
filtrations and is supported on   $K^{2j-a-1}(j,j)$.
\item $\theta_j \from K^{2j-a} \to K^{2j-a+1}$ is a  homeomorphism that respects the filtrations and is supported on   $K^{2j-a}(j,j)$. 
\end{itemize}

Recall the assumption of Proposition 3.1 that $f$ is supported on  $G(b,a)$. By
\cite{Stallings:folding}, $f|G^{\vphantom{\prime}}_a \from G^{\vphantom{\prime}}_a
\to G_a'$ factors as a sequence of folds $\hat \alpha_i$ supported on $G(a,a-1)$ 
followed by a  homeomorphism $\hat \theta_a$ that is simplicial except perhaps on the
$a$-stratum. A fold of edges in  $G_a$ can be viewed as a fold of edges in $G$. Thus
each $\hat \alpha_i$ extends to a fold $\alpha_i$ with support in $G(a,a)$. Define $f_a
\from G \to K^a$ to be the composition of the $\alpha_i$'s. If $\hat \theta_a$ is
simplicial, let $K^{a+1} = K^a$ and let $\theta_a$ be the identity; otherwise $K^{a+1}$
is obtained from $K^a$ by subdividing edges in $K^a(a,a)$ and $\theta_a$ is the
extension of $\hat \theta_a$ over $K^a$. There is an induced homotopy equivalence
$h_{a+1}\from K^{a+1} \to G'$ such that $f = h_{a+1}\theta_af_a$ and such that $h_{a+1}$
is supported on $K^{a+1}(b,a+1)$. Lemma~\ref{push forward} implies that when $K^a$ and
$K^{a+1}$ are equipped with the pushed forward filtrations then $f_a$, $\theta_a$  and
$h_{a+1}$ respect the filtrations. This construction can be repeated, applying it next
to $h_{a+1}$ and continuing by induction, to produce the desired factorization of $f$ as
shown in the following commutative diagram:
\smallskip

Lemma~\ref{homeo} and Lemma~\ref{product}  complete the proof. 
\endproof

We conclude this section with a third application of Lemma~\ref{product}.

\begin{lemma} \label{short} Proposition~\ref{main technical} holds under the additional hypotheses that $f$ is supported on a single stratum $H_j$ and $LC(M_{jj}(f)) < 6$.
\end{lemma}

\begin{remark}  At the end of the proof of Lemma~\ref{j stratum}, it is required that the bound $C = 6$ on  $LC(M_{jj}(f))$ given  above satisfies $C/2 -1 \ge 2$.  This explains our choice of $C = 6$. 
\end{remark}

\proof   A factorization $G = K^0 \xrightarrow{p_1} K^1 \xrightarrow{p_2} \cdots
\xrightarrow{p_k} K^k \xrightarrow{\theta} G'$ of $f$ into folds $p_i$ followed by a
homeomorphism $\theta$ in not unique.  For this proof we factor with a preference for  
$p_i$  that fold into lower strata.  We make this precise as follows.  Let $P_0 := f\from
K^0 \to G'$.  Recall that for any vertex $v \in K^0$ and any oriented edge $e \in
lk(v)$, $DP_0(e) \in lk(P_0(v))$ is the first oriented edge in $P_0(e)$. Edges $e_1, e_2
\in lk(v)$ can be folded if and only if $DP_0(e_1) = DP_0(e_2)$.   If there is a vertex
$v \in K^0_{j-1}$ and an  edge $e_1 \in lk(v) \cap K^0(j,j)$ such that $DP_0(e_1) \in
G'_{j-1}$ then  choose the first fold $p_1 \from K^0 \to K^1$ to use $e_1$ and an edge $e_2
\in K^0 _{j-1}$; this is possible because $P^0|K^0_{j-1} \from K^0_{j-1} \to G'_{j-1}$ is a
simplicial homeomorphism.  Then $P_0 = P_1 p_1$ where  $P_1 \from K^1 \to G'$.  Iterate this
process replacing $P_0$ with $P_1$ and so on to produce  the indicated factorization,
where we denote $\theta p_k \cdots p_{i+1}$ by $P_i \from K^i \to G'$ as shown in the
following commutative diagram:
\smallskip

Let $C_i$ be the sum of all entries of $M_{jj}(P_i)$.  By assumption $C_i$ is uniformly
bounded.  It is helpful to think of $C_i$ as follows.  There is a subdivision of $K^i$
with respect to which $P_i$ is simplicial; we refer to  the edges of this subdivision as
{\it edgelets} of $K^i$.  Edgelets that map into $G'(j,j)$ are colored red while all
other edgelets are colored white.  All the red edgelets are contained in $K^i(j,j)$ and 
$C_i$ is the number of red edgelets.  Two edgelets identified by $p_{i+1}$ must have the
same color.  Thus  $C_{i+1} \le C_i$ for all $i$ and $C_{i} < C_{i-1}$ if  $p_{i+1}$
identifies any red edgelets. In the latter case we say that    {\it $p_{i+1}$ decreases
$j$-length}.   There is a uniform bound to the number of $p_{i+1}$'s that decrease
$j$-length.   

If  $e$ is an edge of $K^{i-1}(j,j)$ then $p_{i}(e)$ is an edge path of length one or two.  If the length is two and if both edges in $p_{i}(e)$ contain red edgelets, then part of $e$ was folded with all of some other edge that contained red edgelets; in particular, some red edgelets of $e$ were identified by $p_{i}$.  In all other cases the edge paths  $p_{i}(e)$ and $e$ have the same number of edges (either zero or one) that contain red edgelets.  We conclude, by induction on $m$, that if $e$ is an edge of $K^{i-1}(j,j)$ and if the edge path $p_mp_{m-1}\cdots p_{i}(e)$ contains two edges  in $K^{m}(j,j)$ that contain red edgelets then  $p_l$ decreases $j$-length for some $i \le l \le m$.  

If $p_{i}$ does not fold into lower strata then $p_{i}(e) \subset K^{i}(j,j)$ for each edge $e \subset K^{i-1}(j,j)$.  There exists $A \ge 0$ with the following property.  If $p_l$ does not fold into lower strata for all $i \le l \le i+A$ then there exists an edge $e$ of $K^{i-1}(j,j)$ such that  the edge path $p_{i+A}p_{i+A-1}\cdots p_{i}(e)$  contains two loops in $K^{i+A}(j,j)$ by 
an argument similar to one used in the proof of Lemma 5.1.  Since every loop in $K^{i+A}(j,j)$ contains red edgelets,  $p_l$ decreases $j$-length for some $i \le l \le i+A$.

    By Lemma~\ref{product}, Lemma~\ref{homeo}, Lemma~\ref{one fold} and Lemma~\ref{pg} it suffices to show that there is a uniform bound to the number of $p_{i}$'s that do not fold into lower strata.  It therefore suffices to show that if  $p_{i}$ does not fold into lower strata then $p_l$ decreases $j$-length for some $i \le l \le i+A$.  

We say that an edge $e$ in $K^{i-1}(j,j)$  with initial vertex $v \in K^{i-1}_{j-1}$ is a {\it frontier edge} and that    $K^{i-1}$ has the {\it frontier red edgelet property} if the initial edgelet of every frontier edge is red.  Assume that $p_{i}$ does not fold into lower strata.  Given our preference for folding into lower strata,  $K^{i-1}$ must have the frontier red edgelet property.  If an initial segment of  a frontier edge is folded by $p_i$ then $p_i$ decreases $j$-length.  If not, then $K^{i}$ has the frontier red edgelet property.     Thus either $p_l$ decreases $j$-length for some $i \le l \le i+A$  or   $p_l$ does not fold into lower strata for all $i \le l \le i+A$.  By our choice of $A$, $p_l$ decreases $j$-length for some $i \le l \le i+A$.
\endproof

\section{Irreducibility}  \label{irreducibility}

A filtration $\hat \F$  of length $N+1$ {\it refines the $j$-stratum} of $\F$ if $\Hat \F$ is obtained from $\F$ by inserting a new filtration element between $G_{j-1}$ and $G_j$.  More precisely, $\hat G_i = G_i$ for $i < j$ and $\hat G_i = G_{i-1}$ for $i >j$.  

Suppose that $f \from G \to G'$ respects the filtrations $\F$ and $\F'$ and restricts to an immersion on each edge of $G$.  We say that  $f$ is {\it $(\F,j)$-reducible} if:
\begin{itemize}
\item  there is a filtration  $\hat \F$ that refines the $j$-stratum of $\F$
\item  $T(f|\hat G_j,v) \ge 2$ for all vertices $v \in \hat G_j$. 
\item  $f|\hat G_j \from \hat G_j \to f(\hat G_j)$ is a homotopy equivalence
\end{itemize}
If  $f$ is not   $(\F,j)$-reducible  then it is  {\it $(\F,j)$-irreducible}.

We will apply $(\F,j)$-irreducibility via the following lemma and corollary.
 
\begin{lemma} \label{homology changes} Suppose that  $f$    is $(\F,j)$-irreducible. Suppose also that   $X$ is a subgraph of $G$ such that   $G_{j-1} \subset X $ and $X \subset  G_j$
are proper inclusions and such that $T(f|X,v) \ge 2$ for all vertices $v \in X$.  Then  either $\rk(H_1(f(X))) > \rk(H_1(X))$ or $\rk(H_0(f(X))) < \rk(H_0(X))$.
\end{lemma}

\proof  The conditions on $X$ imply that $f(X) \subset G_j'$ properly contains $G_{j-1}$ and has no valence one vertices.  If  there is a rank preserving bijection between the components of $X$ and the components of $f(X)$ then $f|X \from X \to f(X)$ is a homotopy equivalence in contradiction to the assumption that $f$ is $(\F,j)$-irreducible.
\endproof

\begin{cor} \label{gamma grows}  Suppose that $K^i$, $0 \le i \le 2n$,  are marked graphs of rank $n$ with filtrations $\F^i$ of length $N$, that $1 \le j \le N$, that  each $f_i \from
K^{i-1} \to K^i$ respects the filtrations and that each $f_i$ is 
$(\F^{i-1},j)$-irreducible.   Suppose  further that $\gamma \subset K^0_j$ is an
immersed loop that is not contained in $K^0_{j-1}$ and that   $f_{2n}\circ \ldots \circ
f_1|\gamma$ is an immersion.  Then $f_{2n}\circ \ldots \circ f_1(\gamma)$ contains
$K^{2n}(j,j)$.
\end{cor}

\proof Let $X_i$ be the smallest subgraph of $K^i$ that contains both  $K^i_{j-1}$ and the immersed loop $f_i\circ \ldots \circ f_1(\gamma)$.  Then $f_{i+1}(X_i) = X_{i+1}$   and if  $X_i \ne K^i_j$, then $X_i$ satisfies the hypotheses of Lemma~\ref{homology changes} with respect to $f_{i+1}\from
K^i \to K^{i+1}$.  The corollary therefore follows from Lemma~\ref{homology changes},
the fact that $\rk(H_1(X_{2n})) \le n$  and the fact that $\rk(H_0(X_0)) \le 
\rk(H_1(X_0))$.   
\endproof

The next lemma is inspired by the following observation.  If $M$   is a product of
positive integer  matrices $M_i$ then $\mlog(LC(M)) \sim \sum \mlog(LC(M_i))$.  This
fails without the assumption that $M$ is positive; for example,  if the $M_i$ are  upper triangular with ones on the diagonal then (c.f.
Lemma~\ref{pg}) $LC(M) \sim \sum LC(M_i)$).

\begin{lemma} \label{j stratum}Suppose that $f\from G \to G'$ restricts to an immersion
on edges, $f$  respects filtrations $\F$ and $\F'$, $f$ is supported on $G(j,j)$, and
$f$  factors as $G = K^0 \xrightarrow{f_1} K^1 \xrightarrow{f_2} \cdots
\xrightarrow{f_k} K^k \xrightarrow{\theta}  G'$  where:
\begin{enumerate}
\item $\theta$ is a   homeomorphism.
\item $f_i$ respects filtrations $\F^{i-1}$, $\F^i$ where $\F^0=\F$ and 
$\theta$ pushes $\F^k$ forward to $\F'$.
\item  each $LC(M_{jj}(f_i)) \ge 6$.
\item   each $f_i\from K^{i-1} \to K^i$  has a homotopy inverse $g_i\from K^i \to K^{i-1}$
satisfying properties (2) and (3) of Proposition~\ref{main technical} and such that
$\log(LC(M_{jj}(f_i))) \succ \mlog (LC(M_{jj}(g_i)))$. 
\item each $f_i$ is $(\F^{i-1},j)$-irreducible.
\end{enumerate}
Let  $\theta' : G' \to K^k$  be the homotopy inverse to $\theta$ given by Lemma~\ref{homeo}.   Then the homotopy inverse $g = g_1 \circ g_2 \circ \cdots g_k \circ \theta' \from G' \to G$  satisfies properties (2) and (3) of
Proposition~\ref{main technical} and 
$$  
\log (LC(M_{jj}(f))) \succ \sum_{i=1}^{k} \log (LC(M_{jj}(f_i))) \succ \mlog (LC(M_{jj}(g)))
$$ 
\end{lemma} 

\proof There is no loss in assuming that $\theta$ is the identity and that $k > 1$.  Items (2) and (3) of
Proposition~\ref{main technical} are closed under composition so are satisfied by $g$. 

 Choose $1 \le m_1 \le k$ so that $LC(M_{jj}(f_{m_1}))$ is maximal.  Choose $m_2$ such that $LC(M_{jj}(f_{m_2}))$ is maximal among all $LC(M_{jj}(f_{i}))$ with $|i-m_1| > 2n$.  Choose $m_3$ such that $LC(M_{jj}(f_{m_3}))$ is maximal among all $LC(M_{jj}(f_{i}))$ with $|i-m_1| > 2n$ and $|i-m_2| > 2n $.  Continue in this way as long as possible to choose $\{m_1,\ldots, m_p\}$.  It is clear that 
$$
\sum_{l=1}^p\log(LC(M_{jj}(f_{m_l}))) < \sum_{i=1}^{k} \log(LC(M_{jj}(f_i))) < (4n+1) \sum_{l=1}^p\log(LC(M_{jj}(f_{m_l})))  
$$
so
$$
\sum_{l=1}^p\log(LC(M_{jj}(f_{m_l}))) \sim \sum_{i=1}^{k} \log(LC(M_{jj}(f_i))). 
$$

Lemma~\ref{direction easy} and item 4 imply that 
$$
\sum_{i=1}^{k} \log(LC(M_{jj}(f_i))) \succ \sum_{i=1}^{k} \mlog (LC(M_{jj}(g_i))) \succ
\log(LC(M_{jj}(g))).  
$$ 

Reorder the $m_l$'s so that they are increasing. Define
$$
\alpha:=  f_{m_1-1} \circ \cdots \circ f_1
$$
$$
\omega:=   f_k \circ \cdots \circ f_{m_p +1}
$$
and for  $1 \le l \le p-1$  
$$
h_l := f_{m_{l+1}-1} \circ \cdots \circ f_{m_l+1}.
$$  
Thus $f$ factors as 
$$
K^0 \xrightarrow{\alpha} K^{m_1-1} \xrightarrow{h_1 \circ f_{m_1} } K^{m_2-1} \xrightarrow{h_2 \circ f_{m_2} } \cdots \xrightarrow{h_{p-1} \circ f_{m_{p-1} } }K^{m_p-1} \xrightarrow{f_{m_p}} K^{m_p} \xrightarrow{\omega}  K^k
$$  

For each $1 \le l \le p-1$, there is an edge $e_{l}$ in  $K^{m_l-1}(j,j)$  and an oriented edge $e_{l}^*$ in $K^{m_l}(j,j)$   that occurs at least  $LC(M_{jj}(f_{m_l}))/2$ times in $f_{m_l}(e_{l})$.  There are $(LC(M_{jj}(f_{m_l}))/2)-1$ subintervals of $e_l$ with disjoint interior whose $f_{m_l}$-images are immersed loops $\{\gamma_u\}$ that contain  $e_{l}^*$.  Since $(h_l \circ f_{m_l}) $ restricts to an immersion on edges, $h_l|\gamma_u$ is an immersion and Corollary~\ref{gamma grows} implies that  each $h_l(\gamma_u)$ contains  $K^{m_l}(j,j)$.  Thus   $(h\circ f_{m_l})(e_{l-1})$ covers every edge of  $K^{m_l}(j,j)$ at least $(LC(M_{jj}(f_{m_l}))/2)-1$ times.  Note that by item 3, $(LC(M_{jj}(f_{m_l}))/2)-1\ge 2$.

It follows  that 
$$
\log (LC(M_{jj}(f))) \succ \sum_{l=1}^p\log(LC(M_{jj}(f_{m_l})))  
$$
which completes the proof.

\endproof

\section{Proof of Proposition~\ref{main technical}}  

    We  now  present the proof of Proposition~\ref{main technical}. Recall that $G$ and $G'$ are marked graphs of rank $n$ with  filtrations $\F$ and $\F'$ of length $N$,  $f :G \to G'$ respects the filtrations, restricts to an immersion on each edge and  is supported on $G(a,b)$.  Moreover,  $G$ has at most $V_n/2$ vertices with $T(f,v) = 2$ and $T(f|G_i,v) \ge 2$ for all vertices $v \in G_i$ and all $i$.  In particular, $G$ has at most $V_n$ vertices and at most $Edge_n$ edges.

The proof of Proposition~\ref{main technical}  is by downward induction on  $N$. Recall that the maximal length of a filtered marked graph is $2n-1$. Proposition~\ref{main technical}  is therefore vacuously satisfied for $N = 2n$ providing the basis step for our downward induction.  We  assume that  $N <  2n$ and  that the proposition holds for filtrations with length greater than $N$.    By Lemma~\ref{one stratum suffices}, we may assume that $f$ is supported on a single stratum $G(j,j)$.

As a special case, suppose that $f$ is  $(\F,j)$-reducible.  
   Choose a refinement $\hat \F$  of $\F$ as in the definition of $(\F,j)$-reducible.  Define $\hat \F'$ by inserting $f(\hat G_j)$ between $G_{j-1}'$ and $G_j'$.  Then $f \from G \to G'$ respects $\hat \F$ and $\hat \F'$ and satisfies the hypotheses of 
Proposition~\ref{main technical}; note that  $f$ is supported on $G(j,j)$ with respect to the original filtrations and is supported on $G(j+1,j)$ with respect to the new filtrations.  By the inductive hypothesis, the conclusions of Proposition~\ref{main technical} hold for the transition matrices defined with respect to   $\hat \F$ and $\hat \F'$.  But for each $r, s$, there exists $r',s'$ such that $M_{rs}(f)$ defined with respect to $\F$ and $\F'$ equals $M_{r's'}(f)$ defined with respect to $\hat \F$ and $\hat \F'$.  This verifies the proposition in the special case.

\begin{sublem} \label{sublemma}  $f$ factors as 
$$G = K^0 \xrightarrow{f_1} K^1 \xrightarrow{f_2} \cdots
\xrightarrow{f_m} K^m \xrightarrow{\theta}  G'
$$   where:
\begin{enumerate}
\item $\theta$ is a   homeomorphism.
\item   Proposition~\ref{main technical} holds  for each $f_i \from K^{i-1} \to K^i$.
\item    $f_i$ is $(\F^i,j)$-irreducible  and $LC(M_{jj}(f_i)) \ge 6$ for $i \le m-1$.
\end{enumerate}
\end{sublem}

\proof     Begin with a factorization $f = \theta \circ p_{w} \circ \cdots \circ p_1$  where each
$p_l$ is a fold and $\theta$ is a homeomorphism. Define a finite sequence of
integers $0=k_0  < k_1  < \ldots < k_m =w$ as follows. If $f$ is $(\F,j)$-reducible or if
$LC(M_{jj}(f)) < 6$, let $k_1 = w$. Otherwise $k_1$ is the smallest positive value for
which $p_{k_1}\circ \ldots\circ p_1 $ is $(\F,j)$-irreducible and for which
$LC(M_{jj}(p_{k_1}\circ \ldots \circ p_1)) \ge 6$. Define $k_2 > k_1$ similarly
replacing $f$ with $\theta \circ p_{w} \circ\ldots p_{k_1+1}$. This process stops in
finite time to produce $0=k_0  < k_1  < \ldots < k_m = w$. Denote $G$ by $K^0$ and
$p_{k_i}\circ \cdots \circ p_{k_{i-1}+1}$ by $f_i\from K^{i-1} \to K^i$. Thus $f =
\theta \circ f_m \circ \cdots \circ f_1$ and it remains to verify the item~2.  

    Suppose at first that $i < m$.  Then  $f_i = p_{k_i} \circ \hat f_i$ where either $\hat f_i$  is $(\F^{i-1},j)$-reducible  or   $LC(M_{jj}(\hat f_i))< 6$.    Lemma~\ref{short} and the special case considered above imply that Proposition~\ref{main technical} holds for $\hat f_i$. Lemma~\ref{one fold} and Lemma~\ref{product} then imply that 
Proposition~\ref{main technical} holds  for  $f_i$ .    If this argument does not also apply to $f_m$, then  either the $f_m$ is $(\F^{m-1},j)$-reducible  or   $LC(M_{jj}(\hat f_m))< 6$.   In either case Proposition~\ref{main technical} holds  for  
$f_m$. 
\endproof  

    Let $g_i \from K^k \to K^{k-1}$ be the   homotopy inverse of $f_i$ determined by item (2) of Sublemma~\ref{sublemma}.   In
particular,
$$
\mlog(LC(M_{rs}(g_i))) \prec \mlog(LC(M_{rs}(f_i))).
$$ 
For all but  uniformly many $i$, $\abs{\Fr(K^{i-1},j,j)} = \abs{\Fr(K^i,j,j)}$.   We
are therefore  reduced  by Lemma~\ref{product} to the case that $\theta$ is the identity,
that  each $g_i$ is supported on $K^k(j,j)$ and that item (3) of Sublemma~\ref{sublemma}  is satisfied
for $i \le m$.  We may assume without loss that $j = s = N$.

Consider the homotopy
commutative diagram
and for $l=1,\ldots,m$ denote homotopy inverses
\begin{align*}
\hat g_l &= g_l \circ \ldots\circ g_m \from K^m \to K^{l-1} \\
\hat f_l &= f_m \circ \ldots \circ f_l \from K^{l-1} \to K^m
\end{align*}
We set notation as follows. For $h \from X \to Y$ supported on $X(j,j)$, let  $A(h)
= M_{j,j}(h)$ and let $B(h)$ be  the submatrix of $M(h)$ whose columns correspond to 
$X(j,j)$ and whose rows correspond to $Y(j-1,r)$.  Note that 
$$
LC(M_{rs}(h)) =\max\{ LC(A(h)), LC(B(h))\}.
$$

Let $g = g_1 \circ g_2 \circ \cdots g_m \from G' \to G$.  Lemma~\ref{j stratum} implies  that $\log(LC(A(g))) \prec \log(LC(A(f)))$ and we are reduced to showing that   $\log(LC(B(g))) \prec \log(LC(M_{rs}(f))$.

Our first estimate follows from Lemma~\ref{direction easy} and  Lemma~\ref{j stratum}.\begin{eqnarray*}
\mlog(LC(A(\hat g_{l-1}))) & \prec &  \sum_{i=l-1}^m \mlog(LC(A(g_i))) \\
                           & \prec &    \sum_{i=l-1}^m \log(LC(A(f_i))) \\
                           & \prec &  \sum_{i=1}^m \log(LC(A(f_i))) \\
                           & \prec &  \log(LC(A(f))) \\
\end{eqnarray*}

For our next estimate, it  is convenient to work without  $\log$.  Keep in mind that  $\log(x) \prec \log(y)$ if and only if $x < y^c$ for some uniform constant $c\ge 1$.
 By definition,
$$
B(g) = B(g_m) + \sum_{l=1}^{m-1}  B(g_l)A(\hat g_{l+1}).
$$
Thus for some uniform constant $c\ge 1$
\begin{eqnarray*}
LC(B(g))  & \le & LC(B(g_m))+ Edge_n  \cdot  \sum_{l=1}^{m-1}  LC(A(\hat g_{l+1})) \cdot (LC(B(g_l))  \\
            & \le & Edge_n  \cdot   LC(A(f))^c \cdot  (\sum_{l=1}^m LC(B(g_l))) \\ 
           & \le &  Edge_n  \cdot  LC(A(f))^c \cdot  (\sum_{l=1}^m LC(M_{rs}(f_l))^c)\\
           & \le &  Edge_n  \cdot  LC(A(f))^c \cdot  (\sum_{l=1}^m LC(M_{rs}(f_l)))^c.  
\end{eqnarray*}
We restate this back in terms of $\log$ and continue, using both  Lemma~\ref{j stratum} and the fact that  $B(f) = B(f_1) + \sum_{l=2}^{m}  B(f_l)A(\hat f_{l-1})$ to verify the next to the last inequality.

\begin{eqnarray*}
\lefteqn{\mlog(LC(B(g)))}\\
  &  \prec  &   \log(LC(A(f))) + \log(\sum_{l=1}^m LC(M_{rs}(f_l))) \\ 
  &  \prec  &   \log(LC(A(f))) + \log(\sum_{l=1}^m LC(A(f_l)) + \sum_{l=1}^m LC(B(f_l))) \\   
  &  \prec  & \log(LC(A(f))) +\sum_{l=1}^m \log(LC(A(f_l))) + 
\log(\sum_{l=1}^m LC(B(f_l))) \\ 
 &  \prec  &  \log(LC(A(f))) + \log(LC(A(f))) + \log(LC(B(f))) \\ 
 &  \prec  &   \log(LC(M_{rs}((f))).
\end{eqnarray*}
which completes the proof.


\begin{thebibliography}{BFH00}

\bibitem[Ali02]{Alibegovic:translation}
Emina Alibegovi{\'c}, \emph{Translation lengths in {${\rm Out}(F\sb n)$}},
  Geom. Dedicata \textbf{92} (2002), 87--93.

\bibitem[BFH00]{BFH:TitsOne}
M.~Bestvina, M.~Feighn, and M.~Handel, \emph{{The Tits alternative for ${\rm
  Out}(F\sb n)$. I. Dynamics of exponentially-growing automorphisms.}}, Ann. of
  Math. \textbf{151} (2000), no.~2, 517--623.

\bibitem[BH92]{BestvinaHandel:tt}
M.~Bestvina and M.~Handel, \emph{Train tracks and automorphisms of free
  groups}, Ann. of Math. \textbf{135} (1992), 1--51.

\bibitem[HM06]{HandelMosher:parageometric}
M.~Handel and L.~Mosher, \emph{Parageometric outer automorphisms of free
  groups}, Trans. AMS (2006), to appear. Preprint, arXiv:math.GR/0410018, 2004.

\bibitem[LL03]{LevittLustig:NorthSouth}
G.~Levitt and M.~Lustig, \emph{Irreducible automorphisms of {$F\sb n$} have
  north-south dynamics on compactified outer space}, J. Inst. Math. Jussieu
  \textbf{2} (2003), no.~1, 59--72.

\bibitem[Sta83]{Stallings:folding}
J.~Stallings, \emph{Topology of finite graphs}, Inv. Math. \textbf{71} (1983),
  551--565.

\bibitem[Vog02]{Vogtmann:OuterSpaceSurvey}
K.~Vogtmann, \emph{Automorphisms of free groups and outer space}, Proceedings
  of the Conference on Geometric and Combinatorial Group Theory, Part I (Haifa,
  2000), vol.~94, 2002, pp.~1--31.

\end{thebibliography}

\providecommand{\bysame}{\leavevmode\hbox to3em{\hrulefill}\thinspace}

\bigskip

\bigskip\noindent
\textsc{\scriptsize
Michael Handel:\\
Department of Mathematics and Computer Science\\
Lehman College - CUNY\\
250 Bedford Park Boulevard W\\
Bronx, NY 10468\\
michael.handel@lehman.cuny.edu
}

\bigskip

\bigskip\noindent
\textsc{\scriptsize
Lee Mosher:\\
Department of Mathematics and Computer Science\\
Rutgers University at Newark\\
Newark, NJ 07102\\
mosher@andromeda.rutgers.edu
}

\end{document}